\documentclass[12pt]{amsart}
\usepackage{amssymb}
\usepackage{amsmath}
\usepackage{latexsym}
\usepackage{amsthm}
\usepackage{graphicx}
 \makeatletter
\def\thmhead@plain#1#2#3{%
 \thmname{#1}\thmnumber{\@ifnotempty{#1}{
 }#2}%
 \thmnote{ \the\thm@notefont(#3)}}
\let\thmhead\thmhead@plain
\def\swappedhead#1#2#3{%
 \thmnumber{#2}\thmname{\@ifnotempty{#2}{. }#1}%
 \thmnote{ \the\thm@notefont(#3)}}
\makeatother

\theoremstyle{definition} %%% for statements in roman typeface
\newtheorem{definition}{Definition}[section]

\newtheorem{example}[definition]{Example}
\newtheorem{examples}[definition]{Examples}

\theoremstyle{plain}      %%% for statements in italic typeface
\newtheorem{proposition}[definition]{Proposition}
\newtheorem{theorem}[definition]{Theorem}

\newtheorem{lemma}[definition]{Lemma}
\newtheorem{fact}[definition]{Fact}
\begin{document}
\address{Clarkson University}
\email{aalhakim@clarkson.edu, akinwamb@clarkson.edu}

\title[De Bruijn Graph Homomorphisms and Recursive De Bruijn Sequences]{De Bruijn Graph Homomorphisms and Recursive De Bruijn Sequences}
\author[Alhakim and Akinwande]{Abbas Alhakim and Mufutau Akinwande \\Department of Mathematics\\ \& Computer Science\\
Clarkson University\\
Potsdam, NY 13699}
% %%%%%%%%%%%%%%%%%%%%%%%%%%%%%%%%%%%%%%%%%%%%%%%%%%%%%%%%%%%%%%%%%%%%%%%%%%%%%%%%%
\begin{abstract}
This paper presents a method to find new De~Bruijn cycles based on ones of lesser order. This is done by mapping a De~Bruijn cycle to several vertex disjoint cycles in a De~Bruijn digraph of higher order and connecting these cycles into one full cycle. We characterize homomorphisms between De~Bruijn digraphs of different orders that allow this construction. These maps generalize the well-known D-morphism of Lempel \cite{Lempel70} between De~Bruijn digraphs of consecutive orders. Also, an efficient recursive algorithm that yields an exponential number of nonbinary De~Bruijn cycles is implemented.
\end{abstract}

\maketitle
\section{Introduction}
The Lempel homomorphism between binary de Bruijn graphs of consecutive orders has been used by many authors for over thirty years to construct new de Bruijn sequences using a given one of smaller order as in \cite{Annex97}, \cite{Chang99}, \cite{Lempel70}, or to obtain results about the linear complexity of binary sequences as, e.g., in \cite{Chan82}, \cite{EtzionLempel84} and references therein.

Ronse gives an attempt in \cite{Ronse84} to generalize Lempel's construction by presenting a single function between two non-binary De~Bruijn digraphs of consecutive orders that has a similar effect. It appears that there is hardly any other attempt in literature to generalize this construction despite the importance and simplicity of the approach. It is worth mentioning that \cite{Chen1995} presents a collection of homomorphisms from higher order De~Bruijn graphs to lower order ones, but these do not nearly enjoy the properties of the Lempel homomorphism as will be shown in the sequel, and thus can not be used to construct De~Bruijn cycles. In this paper we generalize Lempel's homomorphism by describing and characterizing a class of homomorphisms between two De~Bruijn digraphs of arbitrarily different orders but with the same alphabet, the direction of these functions being of course from the higher order digraph to the lower order one.  That is, we both generalize to non-binary alphabets and consider De~Bruijn digraphs of non-consecutive orders. The definition of the proposed class is exceedingly simple. A homomorphism belongs to it if any cycle in the lower order digraph is the image of a fixed number of nonintersecting paths in the higher order digraph and of the same size as the given cycle. This fixed number depends only on the difference between the orders of the digraphs.

Besides the fact that these homomorphisms attest to the structural similarity among different order De~Bruijn digraphs, their significance stems from their applicability to interesting problems. For example, noticing that a generalized feedback shift register sequence (GFSR), described for example in \cite{MK92}, is indeed a cycle in an appropriate De~Bruijn digraph of some large alphabet, the recent work \cite{AA08} exploits the proposed homomorphisms to obtain a virtually endless number of parallel streams of pseudo-random numbers based on a given GFSR. The streams are simply the inverse paths of the GFSR sequence. These streams can be made very much ``uncorrelated'' by selecting a homomorphism that is highly nonlinear and two De~Bruijn digraphs with a large difference between their orders. Furthermore, this is done efficiently by eliminating the traditional preprocessing effort needed for jumping through a fixed stream of pseudo-random numbers.

Another important and relevant problem is the decomposition of a De-Bruijn graph or digraph into a set of disjoint cycles called a factor, see \cite{fred1982}, \cite{Golomb1967}, \cite{Lempel70} for instance.  By definition of the proposed homomorphisms, it can be seen that a factor in the lower order digraph readily induces a factor in the higher order digraph.

In the current work we only deal with the particular case when the vertices of the lower order digraph are all on one cycle, i.e. on a De~Bruijn cycle, the inverse image paths turn out to form a number of non-intersecting cycles that include all the vertices of the higher order digraph. Thus cross-joining these cycles into one big cycle produces a De~Bruijn cycle in the higher order digraph. The focus of this paper is to describe and characterize the proposed class of homomorphisms and to present particular members that simplify the task of cross-joining those inverse images into a new full cycle.

%A De Bruijn sequence of order $n$ is a string of symbols from a finite alphabet such that if you slide a window of size $n$ over it you encounter each of the possible substrings of length $n$ exactly once. These sequences have been used in diverse fields of engineering and computer science, such as the generation of pseudo-random sequences. In fact, linear feedback shift register sequences having the longest period possible are De ~Bruijn sequences without the all zero state. These sequences are called M-sequences or PN-sequences in standard engineering terminology and they have been used in various applications because they are implemented efficiently in hardware devices called linear feedback shift registers (LFSR). Interest in various aspects of these cycles and their generalizations has recently been pushed forward, as they are now seen to belong to a larger family of combinatorial objects, see \cite{Chung92}, \cite{CG04}, \cite{Knuth2005},\cite{Ruskey96} and \cite{Savage97}. We concern ourselves with the classical De Bruijn sequences.

In Section ~\ref{S:prelim} we include definitions and known results that are
useful for the rest of the paper. In Section ~\ref{S:algorithm}, we present an immediate generalization of Lempel's construction to the case of nonbinary alphabet size $q$. This case turns out to be much richer than the binary case as it allows more homomorphisms and more alternating strings to be used for the construction. The number of resulting De~Bruijn cycles of order $n$ is exponential with a huge base of size $q\phi^2(q)$ where $\varphi$ is the Euler totient function.
We describe two variants of an algorithm that performs this cross-join construction. We also present an efficient implementation of this algorithm that is based on a formula that calculates the location of the cross-join.

In Section ~\ref{S:homomorphisms} we introduce the class of homomorphisms discussed above. Proposition ~\ref{P:propD} characterizes all De~Bruijn homomorphisms while Theorem ~\ref{T:genpropertyD} provides a characterization of our proposed class.
Proposition ~\ref{P:span3homomorphism} studies a special homomorphism between two binary De~Bruijn digraphs with an order difference of two, and whose application yields two cycles of different lengths that can be joined relatively easily into one full cycle.

\section{Preliminaries and Basic Results}\label{S:prelim}
\subsection{Terminology}
For positive integers $n$ and $q$ greater than one let $\mathbb{Z}_q^n$ be the set of all $q^n$
vectors of length $n$ with entries in the group $\mathbb{Z}_q$ of residues modulo $q$. When the group structure is not needed we will sometimes refer to elements of this group as symbols. We will interchangeably use the terms
vector, string and word in the sequel to mean a sequence of consecutive symbols. An \textit{order $n$ De ~Bruijn
digraph}, $B_n(q)$, is a directed graph with $\mathbb{Z}_q^n$ as
its vertex set and for two vectors $\textbf{x}=(x_1,\ldots, x_n)$ and $\textbf{y}=(y_1,\ldots, y_n)$, $(\textbf{x};\textbf{y})$ is an edge if and only if $y_i=x_{i+1};\;i=1,\ldots,n-1$. We then say $\textbf{x}$ is a \textit{predecessor} of $\textbf{y}$ and  $\textbf{y}$ is a \textit{successor} of $\textbf{x}$. Evidently, every vertex has exactly $q$ successors and $q$ predecessors. Furthermore, two vertices are \textit{conjugate} if they have the same successors.

A \textit{cycle} in $B_n(q)$ is a path that starts and ends at the same vertex. It is called \textit{vertex disjoint} if it does not cross itself. Two cycles or two paths in the digraph are vertex disjoint if they do not have a common vertex. The \textit{weight} $W(\textbf{w})$ of a word $\textbf{w}$ is the sum--in $\mathbb{Z}_q$--of all elements in $\textbf{w}$. A \textit{translate} of a word (or cycle) $\textbf{c}$ is the word (or cycle) $\textbf{c}+\lambda$ where $\lambda$ is any scalar and addition is component-wise. A cycle is primitive in $B_n(q)$ if it does not simultaneously contain a word $(x_1,\ldots,x_n)$ and any of its translates. For a positive integer $r$ a function $d:\mathbb{Z}_q^r\rightarrow\mathbb{Z}_q$ is said to be \textit{translation invariant} if $d(\textbf{w}+\lambda)=d(\textbf{w})$ for any $\textbf{w}\in\mathbb{Z}_q^r$ and scalar $\lambda$.

A \textit{De~Bruijn cycle} of order $n$ is a Hamiltonian cycle in $B_n(q)$, i.e. a cycle that visits each vertex exactly once. Alternatively, a ``linear'' De ~Bruijn sequence of order $n$ is a sequence of symbols from $\mathbb{Z}_q$ in which every vector of $\mathbb{Z}_q^n$ occurs as a substring exactly
once. For example, $0001110100$ is a binary De~Bruijn sequence of order $3$. Its corresponding De~Bruijn cycle is $000\rightarrow001\rightarrow011\rightarrow111\rightarrow110\rightarrow101\rightarrow010\rightarrow100\rightarrow000$. For the rest of this paper we will denote such a De~Bruijn cycle--and in fact any cycle--as  $[11101000]$,
and consider any rotation of it to be an equivalent cycle. We will denote the linear sequence as $(0001110100)$
and will simply write $B_n$ to denote the binary De Bruijn digraph of order $n$. Finally for the significance and many known algebraic, combinatorial and graph-theoretical methods of construction of De~Bruijn cycles we refer the reader to \cite{fred1982}, \cite{Golomb1967}, \cite{Knuth2005},\cite{Ruskey96}
%to \cite{Chung92}, \cite{CG04}, \cite{fred1982}, \cite{Golomb1967}, \cite{Knuth2005},\cite{Ruskey96}, \cite{Savage97}
and references therein.

\subsection{Lempel's Homomorphism}\label{SS:Lempel} Define a map
$D:\mathbf{B}_n\rightarrow \mathbf{B}_{n-1}$ by
$$D(a_1,\ldots,a_n)=(a_1+a_2,a_2+a_3,\ldots,a_{n-1}+a_n)$$
where addition is modulo $2$. This function defines a graph
homomorphism (see below) and it is known as Lempel's $D$-morphism due to the
fact that it was studied in \cite{Lempel70}, although it
can be traced back to \cite{Leach1960}. Note that
$D(\textbf{x})=D(\textbf{x}+1)$ for all $\textbf{x}\in \mathbb{Z}_2^n$.
We define the dual of a cycle $\mathbf{c}$, to be its bitwise complement $\mathbf{c}+1$. $\mathbf{c}$ is called self-dual if it is a rotation of $\mathbf{c}+1$. The following facts are proved in \cite{Lempel70}.

\begin{fact}\label{LempelFact1} A cycle of length $p$ in $B_{n-1}$ is the $D$-morphic image of two primitive, vertex disjoint cycles of length $p$ in $B_n$ if and only if it has an even number of ones.
\end{fact}
\begin{fact}\label{LempelFact2} A cycle of length $p$ in $B_{n-1}$ is the $D$-morphic image of a self-dual cycle of length $2p$ in $B_n$ if and only if it has an odd number of ones.\end{fact}

Two cycles are called adjacent if a vertex on one cycle has a conjugate on the other cycle. Swapping the successors of these two conjugate words joins the two cycles into one larger cycle. This is why any pair of conjugate words is called a \textit{cross-join pair}, (a similar concept of a cross-join tuple can be defined for $q>2$.)  By Fact \ref{LempelFact1}, if $\mathbf{c}_{n-1}$ is a Hamiltonian cycle in $B_{n-1}$, then every word in $\mathbb{Z}_2^n$ is either on $\mathbf{c}_n$
or on $\mathbf{c}_{n}+1$, the two primitive pre-images of $\mathbf{c}_{n-1}$ by $D$. Lempel
used this idea to construct De ~Bruijn cycles recursively by
rejoining $\mathbf{c}_n$ and $\bar{\mathbf{c}}_{n}$. The most obvious cross-join
pair is the two alternating strings of size $n$ $z_n=010\ldots$ and its
complement $\bar{z}_n$ which can not be on the same cycle.
Recently, this method was implemented in \cite{Annex97}
with an efficient, linear code, and more recently done even more
efficiently in \cite{Chang99} with a jump from a given De Bruijn
cycle in $B_n$ to a higher order $B_{n+k}$, for some integer $k$ that is a power of $2$,
by pre-computing the effect of an iterative application of the
$D$-morphism. %In the rest of this paper we will generalize this technique by presenting a family of homomorphisms from $B_{n+k}(q)$ to $B_{n}(q)$ for any integer $k\geq1$ that can be used similarly to construct De ~Bruijn cycles with higher order recursively ``\`{a} la Lempel''.

%\subsection{The case $k=1$ and $q>2$}\label{S:algorithm}
\section{Generalization To Nonbinary Alphabets}\label{S:algorithm}
The discussion in this paragraph aims at an immediate
generalization to the $q$-ary case of the recursive construction
of full cycles using the Lempel homomorphism as in
\cite{Annex97,Chang99}. In fact, \cite{Ronse84} gives an attempt to generalize this
homomorphism by considering the single function defined by
$\phi(x_1,x_2)=x_2-x_1$. Note that this function is translation
invariant. 

\subsection{Recursive Construction}
We begib with the following proposition which relates translation invariant homomorphisms to primitive cycles. It generalizes Fact~\ref{LempelFact1} and Fact~\ref{LempelFact2} by finding all translation invariant functions that can be used to mimic the Lempel construction of de~Bruijn cycles in the $q$-ary case, as described below.

\begin{proposition}\label{P:recursiveQ}
Let $\Gamma$ be a vertex disjoint cycle of length $l$ in $B_{n-1}(q)$ and $D_{n-1,1}$ be the
homomorphism from $B_{n}(q)$ to $B_{n-1}(q)$ defined as

$$D_{n-1,1}=(d_1(x_1,x_2),d_1(x_2,x_3),\ldots,d_1(x_{n-1},x_n))$$ where
$d_1(x_1,x_2)=\alpha x_1+\beta x_2$, $\alpha$ and $\beta$ are elements in $\mathbb{Z}_q$ such that
$\gcd(\beta,q)=1$ and $\alpha+\beta=0$. Then

(a) $W(\Gamma)=0$ if and only if $\Gamma$ is the image by $D_{n-1,1}$ of a primitive cycle of length $l$ in
$B_{n}(q)$.

(b) $W(\Gamma)=\lambda\neq0$ if and only if $\Gamma$ is the image by $D_{n-1,1}$ of a cycle
\begin{equation}\label{selfdual}
C=\mathbf{c}\cdot(\mathbf{c}+\beta^{-1}\lambda)\cdot(\mathbf{c}+2\beta^{-1}\lambda)\ldots(\mathbf{c}+(r-1)\beta^{-1}\lambda)
\end{equation}
obtained by concatenating $\mathbf{c}$ with its translates, where $r=q/\gcd(\lambda,q)$, $\mathbf{c}=[x_1,\ldots,x_l]$ is an appropriate primitive cycle of length $l$ and $\beta^{-1}$ is the inverse modulo $q$ of $\beta$.
\end{proposition}

\begin{proof}
First, the fact that $D_{n-1,1}$ is a homomorphism is straightforward. The reason for the index notation in $D_{n-1,1}$ and $d_1$ will be clear in the next section.
We will now prove (b). Let $\Gamma=[\gamma_1,\ldots,\gamma_l]$ be an image of a cycle $C$ in $B_{n}(q)$ given by (\ref{selfdual}). Then $\gamma_1+\ldots+\gamma_l=(\alpha x_1+\beta x_2)+\ldots+(\alpha x_{l-1}+\beta x_l)+(\alpha x_l+\beta x_{l+1})= \alpha(x_1+\ldots+x_l)+\beta(x_2+\ldots+x_{l+1})=
(\alpha+\beta)(x_1+\ldots+x_l)+\beta(x_{l+1}-x_1)$. Since $\alpha+\beta=0$ and $x_{l+1}$ is the first element of $(\mathbf{c}+\beta^{-1}\lambda)$ it follows that the last expression is  $\lambda$. Conversely, let $\Gamma=[\gamma_1,\ldots,\gamma_l]$ be a cycle in $B_{n-1}(q)$ with weight $\lambda\neq0$, $x_1$ be an arbitrary value in $\mathbb{Z}_q$ and define $\mathbf{c}$ by letting
\[
x_i=\beta^{-1}(\gamma_1+\ldots+\gamma_{i-1})+x_1,\,i=2,\ldots,l.
\]
Evidently $r$ is the smallest value such that $\mathbf{c}+r\beta^{-1}\lambda=\mathbf{c}$, hence $C$ as defined in (\ref{selfdual}) is an $rl$-cycle in $B_{n+1}(q)$ whose image by $D_{n-1,1}$ is $\Gamma$. To see this note that $\alpha\beta^{-1}=-1$ and for $i=1,\ldots,l-1$, $\alpha x_i+\beta x_{i+1}=\alpha\beta^{-1}(\gamma_1+\ldots+\gamma_{i-1})+\alpha x_1+(\gamma_1+\ldots+\gamma_i)+\beta x_1=\gamma_i$. Also for $j=0,\ldots, r-1$, $\alpha x_{lj+l}+\beta x_{lj+l+1}=\alpha(x_l+j\beta^{-1}\lambda)+\beta(x_1+(j+1)\beta^{-1}\lambda)=\alpha x_l+j\alpha\beta^{-1}\lambda+\beta x_1+(j+1)\lambda=\alpha x_l+\beta x_1+\lambda=
\alpha\beta^{-1}(\gamma_1+\ldots+\gamma_{l-1})+\alpha x_1+\beta x_1+\lambda=-(\gamma_1+\ldots+\gamma_{l-1})+(\alpha+\beta)x_1+\lambda=\gamma_l$. The last equality follows from the definition of $\lambda$ and the assumption $\alpha+\beta=0$. Moreover, note that for $j=0,\ldots,r-~1$, the noncyclic  sequence $(x_{lj+1},\ldots,x_{lj+l}, x_{l(j+1)+1}\ldots, x_{l(j+1)+n})$ does not contain any translates, otherwise two translates would have the same image by $D_{n-1,1}$ violating the assumption that $\Gamma$ is vertex disjoint.

Part (a) can be done similarly by defining $r=1$ and $x_{l+1}=x_1$ in the proof of (b).
%observe that the assumption $\lambda=0$ in the proof of (b) implies that $x_{lj+i}=x_i$ for all $i=1,\ldots, l$ and $ j=0,\ldots,r-1$. That is $C$ is just $r$ rounds of the single cycle $\mathbf{c}$ of length $l$ which is primitive.
\end{proof}

Technically, Part (a) is a special case of (b). We separate them for easy reference below and because Part (b) straightens Proposition ~98(b) in \cite{Ronse84}.  The choice of $\alpha=q-1$ and $\beta=1$ corresponds to Ronse's function. The algorithm below shows how the above proposition helps to construct $q\varphi(q)$  De ~Bruijn cycles of order $n$ with alphabet size $q$ using one De ~Bruijn cycle of order $n-1$ and a homomorphism $D_{n-1,1}$ satisfying the requirements of Proposition~\ref{P:recursiveQ}. A definition and a lemma are in order.
\begin{definition}
For a given value $\lambda\in\mathbb{Z}_q$ we define an alternating sequence $\theta^{(\lambda)}$ to be the infinite sequence $\{e_1,e_2,e_3,\ldots\}$ such that $e_1=0$ and $e_{i+1}=e_i+\lambda$ for $i>1$. Let $\theta^{(\lambda)}_n$ be the suffix of size $n$ of $\theta^{(\lambda)}$. We also define the cycle $\bf{\lambda}_n$ to be the self loop $(\lambda\ldots\lambda;\lambda\ldots\lambda)$ in $B_n(q)$, for any digit $\lambda$.
\end{definition}

Note that the iterates $(i\lambda+\theta^{(\lambda)}_n\rightarrow (i+1)\lambda+\theta^{(\lambda)}_n);\;i=0,\ldots,q-1$ make a cycle of size $q$ in $B_n(q)$ for any $\lambda$ that is coprime with $q$.
\begin{lemma}\label{L:alternating}
For any element $\gamma\neq0$ in $\mathbb{Z}_q$ the inverse image of the self loop $\bf{\gamma}_{n-1}$ (in $B_{n-1}(q)$) by $D_{n-1,1}$ of Proposition ~\ref{P:recursiveQ} is the cycle of length $q$ in $B_n(q)$ formed by the alternating string $\theta^{(\beta^{-1}\gamma)}_n$.
\end{lemma}
\begin{proof}
Let $[0,x_2,\ldots,x_n]$ be  the inverse of $\bf{\gamma}_{n-1}$ that starts with $x_1=0$. We need to show that $x_i=(i-1)\beta^{-1}\gamma$ for $i\geq1$. This is satisfied when $i=1$. Proceeding by induction, assume this is true for $i$. Then $\alpha x_i+\beta x_{i+1}=\gamma$ implies $(i-1)\alpha\beta^{-1}\gamma+\beta x_{i+1}=\gamma$. Since $\alpha\beta^{-1}=-1$ it follows that $x_{i+1}=i\beta^{-1}\gamma$. Finally, since both $\beta$ and $\gamma$ are coprime to $q$ so is the product $\beta^{-1}\gamma$. Thus the cycle formed by $\theta^{(\beta^{-1}\gamma)}_n$ actually has length $q$.
\end{proof}

The following algorithm is a direct generalization of Lempel's binary construction outlined in Subsection~\ref{SS:Lempel}.

\noindent\textbf{Algorithm A}

Given: $a$ and $\lambda\in \mathbb{Z}_q$, ($\lambda$ coprime with $q$), a homomorphism $D_{n-1,1}$ with $\alpha $ and $\beta$ satisfying Proposition ~\ref{P:recursiveQ} and a De ~Bruijn cycle $\Gamma$ of order $n-1$ written linearly as $(\underset{n-1}{\underbrace{\lambda,\ldots,\lambda}},\gamma_n,\ldots\gamma_{q^{n-1}},\underset{n-1}{\underbrace{\lambda,\ldots,\lambda}})$.

\noindent1) Form the inverse image $\{x_{1,j}\},\;j=1,\ldots,q^{n-1}+n$ of $\Gamma$ that starts with $a$.

\noindent2) For $i=2,\ldots,q$ let $\{x_{i,j}\}=\{x_{1,j}\}+(i-1)\beta^{-1}\lambda$.

\noindent3) Regarding the sequences formed above as cycles in $B_{n}(q)$, erase the edge that emerges from the first $n$-vertex in all cycles except the one that starts with $a+(q-1)\beta^{-1}\lambda$.

\noindent4) Erase the last edge (that reaches back to the first string) in all cycles except the starting cycle  $\{x_{1,j}\}$.

\noindent5) For $i=1,\ldots,q-1$, create an edge from the string $(x_{i,1},\ldots,x_{i,n})$ to the one that starts with $(x_{i+1,1},\ldots,x_{i+1,n})$.

\noindent6) For $i=2,\ldots,q$, create an edge from

$(x_{i,q^{n-1}},\ldots,x_{i,q^{n-1}+n-1})$ to $(x_{i-1,2},x_{i-1,3},\ldots,x_{i-1,n+1})$.\vspace{12pt}

There are $q$ inverse images of $\Gamma$ that collectively include each string of size $n$ exactly once. The algorithm outputs one full cycle of order $n$ obtained by connecting the former $q$ cycles. Since $D_{n-1,1}$ is translation invariant, exactly one translate of  $\theta^{(\beta^{-1}\lambda)}_n$  occurs on each separate cycle. By Lemma ~\ref{L:alternating} the first $n$-vertex in each inverse cycle is a translate of $\theta_n^{(\beta^{-1}\lambda)}$. Beginning with $\theta_n^{(\beta^{-1}\lambda)}+a$, Step (5) iterates through all the other translates of $\theta_n^{(\beta^{-1}\lambda)}$, the last one being $\theta_n^{(\beta^{-1}\lambda)}+a+(q-1)\beta^{-1}\lambda$. Since $(x_{i+1,q^{n-1}},\ldots,x_{i+1,q^{n-1}+n-1})$ and $(x_{i,1},\ldots,x_{i,n})$ are both predecessors of $(x_{i+1,1},\ldots,x_{i+1,n})$, it follows that the vertex $(x_{i,2},\ldots,x_{i,n+1})$ is a successor of $(x_{i+1,q^{n-1}},\ldots,x_{i+1,q^{n-1}+n-1})$, this justifies step (6). Steps (3) and (4) of course do the necessary deletions before individual cycles are joined in steps (5) and (6).

Figure~1 illustrates the above construction for $q=5$, $n=2$ and $a=0$ where each row stands for an inverse cycle and the first row represents the cycle in Step (3), the alternating string is viewed as the first vertex in each cycle and the string $(x+i,i,i+1)$ as the last vertex.

\begin{figure}[h]
\includegraphics{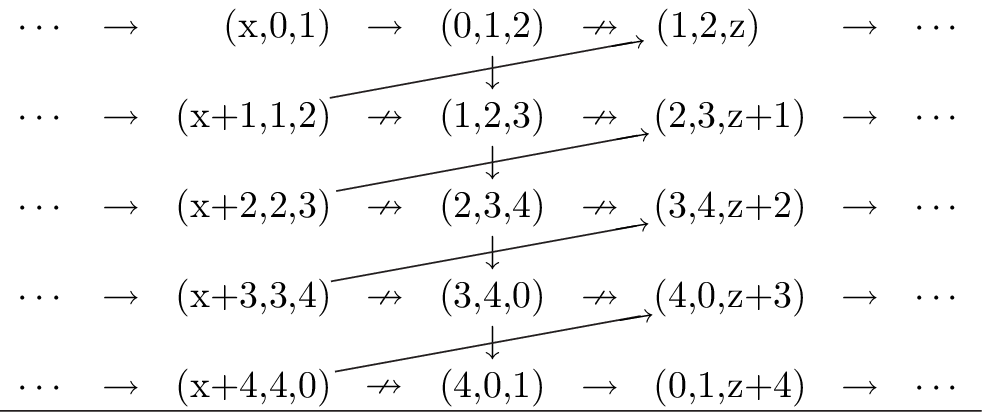}
\caption{}
\end{figure}\label{F:cross-join}

Proposition ~\ref{P:recursiveQ} (b) can also be used to construct the same $q\varphi(q)$ cycles that are the output of Algorithm A. Note that, regardless of the value of $q$, $\Gamma$ has weight zero. Let $\Gamma^{-}$ be the same as $\Gamma$ with exactly one single occurrence of a letter $\gamma$ deleted from the only run of $n-1$ $\gamma$s, where $\gamma$ is chosen such that $\gcd(\gamma, q)=1$. Then the length of $\Gamma^{-}$ as a cycle in $B_{n-1}(q)$ is $q^{n-1}-1$ and $\gcd(W(\Gamma^{-}),q)=1$ and thus constructing the inverse image by $D_{n-1,1}$ of $\Gamma^{-}$ yields a cycle $C^{-}$ in $B_{n}(q)$ of length $q^n-q$. Using Lemma ~\ref{L:alternating} we see that the inverse image of $[\gamma]$ is the cycle formed by iterating $\theta_n^{(\beta^{-1}\lambda)}$.

This shows that the $q$ vertices of $B_n(q)$ that are not in $C^{-}$ form exactly one cycle which is, once cross-joined to $C^{-}$, forms a full cycle. This construction is summarized as follows.\vspace{5pt}

\noindent\textbf{Algorithm B}

\noindent Given: $a$ and $\lambda\in \mathbb{Z}_q$, ($\lambda$ coprime with $q$), a homomorphism $D_{n-1,1}$ with $\alpha $ and $\beta$ satisfying Proposition ~\ref{P:recursiveQ} and a De ~Bruijn cycle $\Gamma$ of order $n-1$ written as $[\underset{n-1}{\underbrace{\lambda,\ldots,\lambda}},\gamma_n,\ldots\gamma_{q^{n-1}}]$.

\noindent1) Let $\Gamma^{-}$ be $[\gamma_n,\ldots,\gamma_{q^{n-1}},\underset{n-2}{\underbrace{\lambda,\ldots,\lambda}}]$ and let $\tilde\Gamma=(\tilde\gamma_1,\ldots,\tilde\gamma_{q^n-q})$ be  a concatenation of $q$ copies of $\Gamma^{-}$.

\noindent2) Let $x_1=a$.

\noindent3) For $i=2,\ldots,n+q-1$ let $x_i=x_{i-1}+\beta^{-1}\lambda$.

\noindent4) For $i=n+q,\ldots,q^n$ choose $x_i$ such that $\alpha x_{i-1}+\beta x_i=\tilde\gamma_{i-n-q+1}$.

Finally, to see that the outcome of this algorithm is the same as that of Algorithm~A the key idea is to note that, for a given $a$ and $\lambda$, the outcomes  of both algorithms start with the $q$ alternating strings in sequence followed by $q$ sections of length $q^{n-1}-1$ each that are all translates of each other.

%%%%%%%%%%%%%%%%%%%%%%%%%%%%%%%%%%%%%%%%%%%%%%%%%%%%%%%%%%%%%%%%%%%%%%%%%%%%%%%
\subsection{Efficient Implementation}
The two algorithms above are mainly theoretical, establishing the possibility of the construction. To actually implement a recursive construction efficiently, note that we need to keep track of the locations of the constant strings within the De~Bruijn sequence being inverted, since the inverse of a constant string is an alternating string where the cross-join must take place. A De~Bruijn sequence that results from input $i$ and $\lambda$ in Algorithm~A will be said to be of type $(i;\lambda)$. Following the notation of Annextein \cite{Annex97}, we say that a sequence representation of a cycle is oriented at $\mathbf{0}_n$ %=\underset{n}{\underbrace{0\ldots 0}}$
if the all zero pattern is written at the end of the sequence, with the understanding that $\mathbf{0}_{n+1}$ is the first vertex in the cycle, e.g., $[11101000]$. With this representation of a cycle $C$, $|\mathbf{x}|_{C}$ denotes the index within $C$ of the ending symbol of $\mathbf{x}$, provided that $\mathbf{x}$ is a substring of $C$.

\begin{proposition}\label{P:Gen_Annex}
Given a homomorphism $D$, let $\Gamma_n$ be any order-n De~Bruijn sequence oriented at $\mathbf{0}_n$, where $n\geq1$. Let $C_i;i=0,\ldots,q-1$ denote the cycle started at symbol $i$ that is the pre-image of $\Gamma_n$ by $D$. Let $\Gamma^{(i;\lambda)}_{n+1}$ be the order-(n+1) De~Bruijn sequence obtained via  Type $(i;\lambda)$. Then for $\gamma\in\mathbb{Z}_q,\gamma\neq0$,

$|\gamma_{n+1}|_{\Gamma_{n+1}^{(i;\lambda)}}=
\begin{cases}
        (q - m)q^n + m                      &   i=0\\
        (m^{\prime}-m)(q^n - 1)             &   i\neq0 \textup{ and } m< m^{\prime} \\
        (m^{\prime} - m)(q^n - 1) + q^{n+1} &   i\neq0 \textup{ and } m > m^{\prime}\\
\end{cases}$

where $m^{\prime} = (- i)\lambda^{-1}$ and $m = (\gamma - i)\lambda^{-1}$, all calculations being in $\mathbb{Z}_q$.
\end{proposition}
In the proof below,  phrases like higher cycle, next cycle, and top cycle refer to the diagrams in Figure~1 and Table~1.
\begin{proof}
Let the homomorphism be induced by the function $d(x_1,x_2)=\alpha x_1+\beta x_2$. The constant $n$-string $\beta\lambda$, denoted $(\beta\lambda)_n$, has an inverse image in each cycle $C_i$. We denote by $^{i}\theta_{n+1}^{(\lambda)}$ the one that belongs to $C_i$. To prove
(i) note that there is an element $m$ such that $\gamma=m\lambda$ for some $m$, since $\lambda$ is relatively prime to $q$. Thus we have $m = \gamma{\lambda}^{-1}$ and $\gamma_{n+1}$ resides on the cycle $C_{m\lambda}$. Furthermore, it is the last string of size $n+1$ on this cycle because $C_0$ is oriented at $\mathbf{0}_{n+1}$.

\begin{table}[h]
%\begin{minipage}[c]{2000mm}
\begin{tabular}{l|c|r}
\hline
\begin{tabular}{ccc}
               ${\longrightarrow}$  & $C_0$              & ${\longrightarrow}$\\
                                    & ${\downarrow}$\\
               ${\longrightarrow}$  & $C_{\lambda}$      & ${\longrightarrow}$\\
                                    & ${\downarrow}$\\
               ${\longrightarrow}$  & $C_{m \lambda}$    & ${\longrightarrow}$\\
                                    & ${\downarrow}$\\
                                    & $\vdots$\\
                                    & (a)\\
               %${\longrightarrow}$  & $C_{(q-1)\lambda}$ & ${\longrightarrow}$\\
\end{tabular}
&
\begin{tabular}{ccc}
               ${\longrightarrow}$  & $C_i$               & ${\longrightarrow}$\\
                                    & $\vdots$\\
               ${\longrightarrow}$  & $C_{0}$      & ${\longrightarrow}$\\
                                    & $\vdots$\\
               ${\longrightarrow}$  & $C_{\gamma}$     & ${\longrightarrow}$\\
                                    & ${\downarrow}$\\
                                    & $\vdots$\\
                                    & (b)\\
\end{tabular}
&
\begin{tabular}{ccc}
               ${\longrightarrow}$  & $C_i$               & ${\longrightarrow}$\\
                                    & $\vdots$\\
               ${\longrightarrow}$  & $C_{\gamma}$      & ${\longrightarrow}$\\
                                    & $\vdots$\\
               ${\longrightarrow}$  & $C_{0}$     & ${\longrightarrow}$\\
                                    & ${\downarrow}$\\
                                    & $\vdots$\\
                                    & (c)\\
\end{tabular}\\
\hline
\end{tabular}\vspace{3pt}
\caption{(a) Type $(0;\lambda)$. (b) Type $(i;\lambda)$ with $l_{\lambda}(i,0) < l_{\lambda}(i,\gamma)$. (c) Type $(i;\lambda)$ with $l_{\lambda}(i,0) > l_{\lambda}(i,\gamma)$.}
\end{table}

\noindent Now, Type $(0,\lambda)$ uses the first $|^0\theta_{n+1}^{(\lambda)}|_{C_0}$ and then proceeds through the alternating strings by adding one digit from each of the cycles $C_\lambda, C_{2\lambda}, ... , C_{m\lambda}$ ($m$ digits) and down to the last cycle.

\noindent All the $(q-m-1)$ cycles that are visited after $C_{m\lambda}$ is first traversed are used in full, thus adding $(q-m-1)q^n$ digits to the length of $|\gamma_{n+1}|_{\Gamma_{n+1}^{(0;\lambda)}}$.

\noindent Finally, the second and final visit to $C_{m\lambda}$ happens just after the location of the alternating string $\theta_{n+1}^{(\lambda)}$. Since $|^0\theta_{n+1}^{(\lambda)}|_{C_{m\lambda}}=|^0\theta_{n+1}^{(\lambda)}|_{C_0}$, we use
$q^n-|^0\theta_{n+1}^{(\lambda)}|_{C_0}$ digits of $C_{m\lambda}$ after this visit.
By adding all the quantities together, we get

$|\gamma_{n+1}|_{\Gamma_{n+1}^{(0;\lambda)}} = |^0\theta_{n+1}^{\lambda}|_{C_0}+m +(q-m-1)q^n+q^n-|^0\theta^{(\lambda)}_{n+1}|_{C_0} = m+(q-m)q^n$.

(ii) When $i\neq0$, type $(i;\lambda)$ is obtained by starting with the cycle $C_i$, as in Table~1 (b) and (c), and connecting through the alternating strings $\theta^{(\lambda)}_{n+1}$. By definition, $|\gamma_{n+1}|_{\Gamma_{n+1}^{(0;\lambda)}}$ is the distance from the start in a cycle oriented at $\mathbf{0}_{n+1}$. For $j\in\mathbb{Z}_q$, let the number of cycles that separate $C_j$ from $C_i$ through steps of size $\lambda$ be denoted $l_{\lambda}(i,j)$. It is easy to see that $l_{\lambda}(i,j)=(j-i)\lambda^{-1}$. We will consider two cases. First, assume $l_{\lambda}(i,0)< l_{\lambda}(i,\gamma)$. Then since $C_i$ is the top cycle, the start of $C_0$ is encountered during the upward traversal. Starting the count at the beginning of $C_0$, we scan $|^0\theta_{n+1}^{(\lambda)}|_{C_0}-1$ digits, i.e., just before the end of the alternating string. The next set of digits are from $C_{-\lambda}$, starting just after the end of the alternating string till the end of this cycle, thus using $(q^n-|^{-\lambda}\theta_{n+1}^{(\lambda)}|_{C_{-\lambda}})$ digits. To reach the constant string $\gamma_{n+1}$ we need to repeat the previous step for all $l_{\lambda}(i,0)$ cycles above $C_0$. Once this is done, the start of $C_i$ is reached. So then we reach the end of $^i\theta^{(\lambda)}_{n+1}$, with $|^i\theta_{n+1}^{(\lambda)}|_{C_i}$ digits and go down along the last digit of each alternating string, down to cycle $C_{i+(q-1)\lambda}$, this adds $(q-1)$ digits. Next, $\gamma_{n+1}$ is located at the end of cycle $C_{\gamma}$ so we need to scan upwards the $(q-l_{\lambda}(i,\gamma)-1)$ cycles underneath $C_{\gamma}$, using $(q^n-1)$ digits of each. Finally, we scan the last  $(q^n-|^{\gamma}\theta_{n+1}^{(\lambda)}|_{C_{\gamma}})$ digits of $C_{\gamma}$ to reach $\gamma_{n+1}$.

Hence, it follows that $|\gamma_{n+1}|_{\Gamma_{n+1}^{(i;\lambda)}}$ is the sum

\[
\left(l_{\lambda}(i,0)(q^n-1)\right)+\left(|^i\theta_{n+1}^{(\lambda)}|_{C_i}\right)+\left(q-1\right)
+\left((q-l_{\lambda}(i,\gamma)-1)(q^n-1)\right)+\left(q^n-|^{\gamma}\theta_{n+1}^{(\lambda)}|_{C_{\gamma}}\right)
\]

Using the fact that $|^j\theta_{n+1}^{(\lambda)}|_{C_j} = |(\beta\lambda)_n|_{\Gamma_n}$ for all $j$, this simplifies to

$(l_{\lambda}(i,0)-l_{\lambda}(i,\gamma))(q^n-1)+q^{n+1}$.

The second case is $l_{\lambda}(i,0)>l_{\lambda}(i,\gamma)$. Then  $C_{\gamma}$ is above $C_0$, see Table~1 (c). Starting at the beginning of $C_0$ we reach the end of $C_{\gamma}$ by scanning upwards $(q^n-1)$ digits of the $(l_{\lambda}(i,0)-l_{\lambda}(i,\gamma))$ cycles. This completes the proof.
\end{proof}

The most important feature of Proposition~\ref{P:Gen_Annex} is that the distance $|\gamma_{n+1}|_{\Gamma_{n+1}^{(i;\lambda)}}$ depends only on the parameters $\gamma$, $\lambda$ and $i$ but not on any previously constructed cycles. The homomorphism parameter $\beta$ is only needed to determine which constant string in $\Gamma_{n+1}$ is the image of $\theta^{(\lambda)}_{n+1}$.

Now we can write down a pseudo-code that is capable of constructing an exponential number of order $n$ De~Bruijn cycles for any alphabet size  $q$, based on the trivial low order cycle. We will use the order-$1$ cycle $\Gamma_1=[12\ldots(q-1)0]$ for odd $q$. Note that the weight $W(\Gamma_1)=0$ if and only if $q$ is odd, thus allowing the construction in algorithm A. Note that the intermediate cycles can be `lifted' to higher orders via distinct homomorphisms. Hence the parameter $\beta$ can be changed in $\phi(q)$ ways. Since there are $q\phi(q)$ type $(i;\lambda)$ cycles for a given homomorphism, we can construct $(q\phi(q)^2)^{n-1}$ cycles of order $n$ using $\Gamma_1$. The following algorithm is an implementation of Algorithm~A above. It takes as input an alphabet size $q$ and three $q$-ary strings $B=(\beta_2,\ldots,\beta_n)$, $L=(\lambda_2,\ldots,\lambda_n)$, $I=(i_2,\ldots,i_n)$. For each recursive step $j$, $\beta_j$ determines the homomorphism to be used for this step, while $(i_j;\lambda_j)$ determines the type of cross-join.

\noindent\textbf{Algorithm AA}

\verb"Input": $q$, and three strings $B=(\beta_2,\ldots,\beta_n)$, $L=(\lambda_2,\ldots,\lambda_n)$, $I=(i_2,\ldots,i_n)$

\verb"Output": a unique order-$n$ De~Bruijn cycle $\Gamma_n$

(1) Let $i_1=\lambda_1=1$ and $\Gamma_1^{(1;1)}=12\ldots(q-1)0$

(2) \verb"For" $j=2,\ldots,n$ repeat the following steps:

(3) Let $\alpha=q-\beta_j$ and calculate $C_0$ as the inverse of $\Gamma_{j-1}^{(i_{j-1};\lambda_{j-1})}$
by the homomorphism induced by $d(x_1,x_2)=\alpha x_1+\beta_j x_2$

(4) Construct $C_k=C_0+k;k=1,\ldots,q-1$, [component-wise addition by $k$]

(5) Let $\gamma=\beta_{j}\times\lambda_{j} (\mod q)$

[Steps (6)-(10) determine the location \verb"pos" of the cross-join within each cycle $C_k$]

(6)  If $j=2$ then let \verb"pos" $=\gamma$ and skip to Step (11)

(7) Let  $m=(\gamma-i_{j-1})\lambda^{-1}_{j-1} (\mod q)$

(8) If $i_{j-1}=0$ let \verb"pos" $=(q-m)q^{j-2}+m$ and skip to Step (11)

(9) Let $m^{\prime}=(-i_{j-1})\lambda^{-1}_{j-1} (\mod q)$

(10) If $m<m^{\prime}$ let \verb"pos" $=(m^{\prime}-m)(q^{j-2}-1)$

\indent otherwise \verb"pos" $=(m^{\prime}-m)(q^{j-2}-1)+q^{j-1}$

[Steps (11) to (23) connect $C_0,\ldots,C_{q-1}$ by the string $\theta^{\lambda}_{j}$]

(11) Let \verb"newstring" be an empty string

(12) If $i_j=0$ do steps (13) to (16) otherwise do steps (17) to (22)

(13) Append the prefix of size \verb"pos" to \verb"newstring"

(14) For $k=1$ to $(q-2)$ append the \verb"pos"$^{th}$ digit of the cycle $C_{k\lambda_j}$ to \verb"newstring"

(15) Append the last $(q^{j-1}-$ \verb"pos"$+1)$ digits of $C_{(q-1)\lambda_j}$ to \verb"newstring"

(16) For $k=(q-1)$ down to $1$

\indent\indent Append the first (\verb"pos" -1) digits of $C_{k\lambda_j}$ to \verb"newstring"

\indent\indent Append the last ($q^{j-1}-$pos) digits of $C_{(k-1)\lambda_j}$ to \verb"newstring"

(17) Let $m^{\prime}=-i_j\lambda_j^{-1} (\mod q)$

(18) For $k=0$ to $m^{\prime}-1$

\indent\indent Append the first (\verb"pos" -1) digits of $C_{-k\lambda_j}$ to \verb"newstring"

\indent\indent Append the last ($q^{j-1}-$pos) digits of $C_{(-k-1)\lambda_j}$ to \verb"newstring"

(19) Append the first \verb"pos" digits of $C_{i_j}$ to \verb"newstring"

(20) For $k=1$ to $(q-2)$ append the \verb"pos"$^{th}$ digit of the cycle $C_{i_j+k\lambda_j}$ to \verb"newstring"

(21) Append the last ($q^{j-1}-$\verb"pos"+1) digits of $C_{i_j+(q-1)\lambda_j}$ to \verb"newstring"

(22) For $k=(q-1)$ down to $(q-1-m^{\prime})$

\indent\indent Append the first (\verb"pos" -1) digits of $C_{k\lambda_j}$ to \verb"newstring"

\indent\indent Append the last ($q^{j-1}-$ \verb"pos") digits of $C_{(k-1)\lambda_j}$ to \verb"newstring"

(23) $\Gamma_j^{(i_j;\lambda_j)}=$ \verb"newstring".

An implementation of this algorithm in C++ is available online, see \cite{AA08w}.  We chose to use odd $q$ because the position, \verb"pos", of the cross-join becomes especially simple to evaluate, see step (6). For even $q$, only the base case $\Gamma_1$ has to be changed to some order-two cycle $\Gamma_2$ and the position of the constant string $\bf{\gamma_2}$ has to be located by searching the string $\Gamma_2$. Table~2 shows all possible output that can be constructed from $\Gamma_1=[120]$ for $q=3$.

\begin{table}[h!]
\small
\begin{tabular}{l|c||l|c}
\hline
$(\beta,i,\lambda)$  &  $\Gamma_2^{(i;\lambda)}$   &   $(\beta,i,\lambda)$  &  $\Gamma_2^{(i;\lambda)}$\\\hline
(1,0,1)              &  [120221100]                &   (2,0,1)              &  [201221100]\\
(1,0,2)              &  [102112200]                &   (2,0,2)              &  [210112200]\\
(1,1,1)              &  [221120100]                &   (2,1,1)              &  [221101200]\\
(1,1,2)              &  [112102200]                &   (2,1,2)              &  [110212200]\\
(1,2,1)              &  [220121100]                &   (2,2,1)              &  [221201100]\\
(1,2,2)              &  [112202100]                &   (2,2,2)              &  [112210200]\\\hline
\end{tabular}
\caption{}
\end{table}

\section{De Bruijn Graph Homomorphisms}\label{S:homomorphisms}
In this section, we study homomorphisms between De ~Bruijn digraphs of
different orders. In Subsection ~\ref{S:propertyD}  we introduce and characterize homomorphisms between
$B_{n+k}(q)$ and $B_n(q)$ for any integer $k\geq1$ that perform like the homomorphisms of the previous section in the sense that taking the inverse by one such homomorphism of a De~Bruijn cycle in $B_n(q)$ dissects $B_{n+k}(q)$ into adjacent vertex disjoint cycles. In Subsection ~\ref{S:binaryPropD}
we single out the binary case which, due to its simplicity, admits a more concise characterization.

\subsection{Homomorphisms with Property (D)}\label{S:propertyD}
Recall that a graph homomorphism is a mapping that preserves the graph structure, so that an edge in $B_{n+k}(q)$ is mapped to an edge in $B_n(q)$. The following proposition characterizes graph homomorphisms between De ~Bruijn digraphs of different orders.

\begin{proposition}\label{P:propD} A necessary and sufficient condition for a
map $D_{n,k}:B_{n+k}(q)\rightarrow B_n(q)$ to be a graph homomorphism is that
\[
D_{n,k}(\mathbf{x})=\left(d_k(x_1,\ldots,x_{k+1}),d_k(x_2,\ldots,x_{k+2}),\ldots,d_k(x_n,\ldots,x_{n+k})\right)
\]
where $\mathbf{x}=(x_1,\ldots,x_{n+k})$ and $d_k$ is any fixed function of $k+1$ variables.
\end{proposition}
\begin{proof}
Sufficiency is quite obvious so we will only prove the necessity. Let $D_{n,k}(x_1,\ldots,\\x_{n+k})$ be $(\tilde{x}_1,\ldots,\tilde{x}_n)$ where
$\tilde{x}_i=d_i(x_1,\ldots,x_{n+k})$
is a function from $\mathbb{Z}_q^n$ to $\mathbb{Z}_q$ for all $i=1,\ldots,n$.
For all values of $x_1,\ldots, x_{n+k},x_{n+k+1}$ $(x_1,\ldots,x_{n+k})$ is a predecessor of $(x_2,\ldots,x_{n+k+1})$. Hence,  since $D_{n,k}$ is a graph homomorphism the diagram below commutes, where the horizontal arrows indicate an edge in the De~Bruijn digraph. That is,
\begin{table}[h]
\center
\begin{tabular}{ccc}
$(x_1,\ldots,x_{n+k})$ & $\overset{B_{n+k}}{\longrightarrow}$ & $(x_2,\ldots,x_{n+k},x_{n+k+1})$\\
$D_{n,k}\downarrow$           &               & $\;\;\downarrow D_{n,k}$\\
$(\tilde{x}_1,\ldots,\tilde{x}_n)$ & $\underset{B_n}{\longrightarrow}$ & $(\tilde{x}_2,\ldots,\tilde{x}_n,\tilde{x}_{n+1})$
\end{tabular}
\end{table}
\begin{equation}\label{E:iterativeEq d_i}
d_i(x_1,\ldots,x_{n+k})=d_{i-1}(x_2,\ldots,x_{n+k+1}); i=2,\ldots,n.
\end{equation}
To finish the proof we need to establish that (1) $d_i(x_1,\ldots,x_{n+k})=d_j(x_1,\ldots,x_{n+k})$ for all $i\neq j$ and (2) $d_i$ depends at most on $x_i,\ldots, x_{i+k}$. We establish this by iterating Equations~ (\ref{E:iterativeEq d_i}) for $i=2,\ldots,n$. To avoid confusion we will denote $d_i$ by $d_i^L$ and $d_i^R$ when it is applied to $(x_1,\ldots,x_{n+k})$ and $(x_2,\ldots,x_{n+k+1})$ respectively (the left and right sides of the diagram above). This is meant to remind us that, e.g., the first variable of $d_i^R$ is $x_2$.

First $\tilde{x}_2=d_2^L(x_1,\ldots,x_{n+k})=d_1^R(x_2,\ldots,x_{n+k+1})$, so that $d_1$ does not depend on its $(n+k)^{th}$ variable and $d_2$ does not depend on its first variable. Next, $\tilde{x}_3=d_3^L(x_1,\ldots,x_{n+k})=d_2^R(x_3,\ldots,x_{n+k+1})$ (noting that by the previous result $d_2^R$ does not depend on its first variable $x_2$). It follows that $d_2$ does not depend on its $(n+k)^{th}$ variable and $d_3$ does not depend on its first and second variables.

Continuing this way we see that for $i=2,\ldots,n$, $d_{i-1}$ does not depend on the last variable and $d_i$ does not depend on its first $(i-1)$ variables. In particular, $d_n$ depends on at most $x_n,\ldots,x_{n+k}$ and $d_{n-1}$ depends on at most $x_{n-1},\ldots,x_{n+k-1}$. Next $d_{n-1}^L(x_{n-1},\ldots,x_{n+k-1})=d_{n-2}^R(x_{n-1},\ldots,x_{n+k})$ implies that $d_{n-2}$ does not depend on its $(n+k-1)^{st}$ variable. Continuing with Equations ~(\ref{E:iterativeEq d_i}) iteratively and backwards this time we establish requirement (2) above. But then Equations ~\ref{E:iterativeEq d_i} reads $d_i^L(x_i,\ldots,x_{i+k})=d_{i-1}^R(x_i,\ldots,x_{i+k})$. Hence, for all $i$ $d_i$ is a fixed function of $k+1$ variables which establishes (1).
\end{proof}

By Fact \ref{LempelFact1} a vertex disjoint cycle in $B_n$ is the $D$-morphic
image of two vertex disjoint cycles in $B_{n+1}$ starting
respectively with zero and one. We generalize this $D$-morphism by finding and characterizing
homomorphisms $D_{n,k}$ so that each vertex disjoint cycle in $B_n(q)$ is the image of exactly
$q^k$ vertex disjoint paths in $B_{n+k}(q)$, one for each
starting string of size $k$. Such a homomorphism (and the
corresponding function $d_k$) \textit{will be said to have
property (D)}.
We will illustrate this property with some examples before we state Theorem
~\ref{T:genpropertyD} which characterizes functions $d_k$
that have property (D). In fact a direct inspection of the
sixteen Boolean functions of two variables shows that the only
homomorphisms with property (D) from $B_{n+1}$ to $B_n$ are the
D-morphism and its bitwise complement.

\begin{examples} (a) Consider the mapping $D_{1,2}$ from $B_{3}$ to $B_1$ that uses the function $d(x_1,x_2,x_3)=x_1+x_2$. The inverse sets of $0$ and $1$ are respectively
$\{000,001,110,111\}$ and $\{010,011,100,101\}$. The edge $(0,1)$
of $B_1$ is mapped back to the four edges
\[\{(001;010),(001;011),(110;100),(110;101)\}\] Note that, even though each edge in $B_1$ is the image of four edges in $B_3$, it is not possible to construct an edge starting with arbitrary strings of size two that is mapped to a given edge of $B_1$. For instance, there is no edge in $B_3$ that starts with either $01$ or $10$ and whose image is the edge $(0;1)$.

\noindent (b) The function $H_{n,k}$ from $B_{n+k}(q)$ to $B_n(q)$ for $k\geq0$ and $n\geq1$ was defined in \cite{Chen1995} as $H_{n,k}(x_1,\cdots,x_{n+k})=(x_{k+1},\cdots,x_{n+k})$. In other words, this function trims the $k$ leftmost symbols of a word so as to make it a word of size $n$. Obviously, this is a homomorphism having, according to the notation of Proposition~\ref{P:propD}, $d(x_1,\cdots,x_{k+1})=x_{k+1}$, hence Theorem~\ref{T:genpropertyD} below shows that it does not enjoy property~(D). In fact the $q^k$ inverses of any cycle in $B_n(q)$ by $H_{n,k}$ disagree only in their first $k$ terms while the body of the sequences are all \textit{equal} to the original cycle. %The reader is invited to find the four inverse sequences of $\mathbf{b}_4$ in Example~\ref{E:toy} by $H_{2,2}$.

\noindent (c) Using the function $d^{(1)}(x_1,x_2,x_3)=x_1+x_3$ however, the edges
$(0;0)$, $(0;1)$, $(1;0)$, $(1;1)$ of $B_1$ are respectively mapped back to
the following sets whose union constitutes the edge set of
$B_{3}$, each edge appearing exactly once.
\begin{eqnarray*}
\{(000;000),(010;101),(101;010),(111;111)\}\\
\{(000;001),(010;100),(101;011),(111;110)\}\\
\{(001;010),(011;111),(100;000),(110;101)\}\\
\{(001;011),(011;110),(100;001),(110;100)\}
\end{eqnarray*}
Hence $d^{(1)}$ enjoys property (D) while $d$ does not. %Table \ref{T:dk} lists all four functions $d^{(1)},\ldots,d^{(4)}$ with property (D) for $k=2$.
\end{examples}
%\begin{table}[h]\label{T:dk}
%\center
%\begin{tabular}{lll|lllll}
%$x_1$ & $x_2$ & $x_3$ & $d$ & $d^{(1)}$ & $d^{(2)}$ & $d^{(3)}$ & $d^{(4)}$\\
%\hline
%0 & 0 & 0 & 0 & 0 & 0 & 1 & 1\\
%0 & 0 & 1 & 0 & 1 & 1 & 0 & 0\\
%0 & 1 & 0 & 1 & 0 & 1 & 1 & 0\\
%0 & 1 & 1 & 1 & 1 & 0 & 0 & 1\\
%1 & 0 & 0 & 1 & 1 & 1 & 0 & 0\\
%1 & 0 & 1 & 1 & 0 & 0 & 1 & 1\\
%1 & 1 & 0 & 0 & 1 & 0 & 0 & 1\\
%1 & 1 & 1 & 0 & 0 & 1 & 1 & 0\\
%\hline
%\end{tabular}\caption{}
%\end{table}
\begin{theorem}\label{T:genpropertyD}
(a) A homomorphism from $B_{n+k}(q)$ to $B_{n}(q)$ that is induced
by $d_{k}(x_1,\ldots,\\x_{k+1})$ enjoys property (D) if and
only if $d_{k}$ is one to one in each of the variables $x_1$ and $x_{k+1}$ when all the other
variables are kept fixed. That is, if and only if $d_{k}(x_1,\ldots,x_{k+1})$ defines a Latin square for each set of fixed values of $x_2,\cdots,x_k$.

(b) The total number of homomorphisms with property (D) is
$(A_q)^{q^{k-2}}$ where $A_q$ is the number of $q\times q$ Latin
squares.
\end{theorem}

\begin{proof}
Part (b) follows directly from (a), see \cite{Sloane} for more about
the sequence $A_q$.  To prove Part (a), first let $d_k$ be a
function with property (D) and $C=[c_1,\ldots,c_l]$ be an arbitrary vertex
disjoint cycle in $B_n(q)$. By definition of property (D) each word $(x_1,\ldots,x_k)$ in $\mathbb{Z}_q^k$
can be appended by a symbol $x_{k+1}$ so that $d_k(x_1,\ldots,x_k,x_{k+1})=c_1$. This says that $d_k$ is surjective from $\mathbb{Z}_q$ to $\mathbb{Z}_q$ (hence injective) with respect to the last variable.

Now let $x_1^{\prime}$ be such that $d_k(x_1^{\prime},x_2,\ldots,x_k,x_{k+1})$ and $d_k(x_1,x_2,\ldots,x_k,x_{k+1})$ are equal to $c_1$.
Since $d_k$ is bijective with respect to the last variable, there exist unique values $x_{k+2},\ldots,x_{n+k+1}$ such that $D_{n,k}(x_2,\ldots,x_{n+k+1})$ $=(c_2,\ldots,c_{n+1})$. If $x_1^{\prime}\neq x_1$ then the two distinct inverse edges \[(x_1,\ldots,x_{n+k};x_2,\ldots,x_{n+k+1}), (x_1^{\prime},x_2\ldots,x_{n+k};x_2,\ldots,x_{n+k+1})\] share a common vertex, contradicting property (D). Hence $d_k$ is one-to-one in the first variable. This establishes the necessary condition.

Conversely, let $d_k$ have the claimed form and let
$\mathbf{c}=[c_1,\ldots,c_l]$ be a vertex disjoint cycle in $B_n(q)$. Given any
string $x_1,\ldots,x_k$ it is possible to find a value
$b\in\mathbb{Z}_q$ so that $d_k(x_1,\ldots,x_k,b)=c_1$, since $d_k$ is
surjective with respect to the last variable. Hence the value $c_1$ has
a set of $q^k$ inverse images that includes all possible strings of size $k$ as
prefixes. The same argument implies that $\mathbf{c}$ has exactly $q^k$ inverse images. To show property
(D) we need to show that no substring of size $n+k$ occurs more
than once in the collection of pre-images. Write the pre-images of
$\mathbf{c}$ as a rectangular array $(x_{ij}); i=1,\ldots,q^k,j=1,\ldots,
k+l+n-1$ (corresponding to the ``linearized'' cycle $c_1,\ldots,c_l,c_1,\ldots,c_{n-1}$)
where the set of prefixes of size $k$ coincides with the
$q^k$ distinct words of this size. Let us denote by $\omega_{ij}(u)$ the
substring of size $u$ on the $i^{th}$ row of $(x_{ij})$ that
starts with $x_{ij}$. Assume there exist integers
$i_1,i_2,j_1,j_2$ such that $\omega_{i_1,j_1}(n+k)$ coincides with
$\omega_{i_2,j_2}(n+k)$. Obviously $j_1\neq j_2$ implies that a
string of size $n$ occurs twice in $\mathbf{c}$, thus contradicting the
assumption that it is vertex disjoint. Assume then that
$j_0$ is the smallest integer with
$\omega_{i_1,j_0}(n+k)=\omega_{i_2,j_0}(n+k)$, which in particular
means that $\omega_{i_1,j_0}(k)=\omega_{i_2,j_0}(k)$. By
construction of $(x_{ij}),\,\,j_0>1$. Hence
$d_k(\omega_{i_1,j_0-1}(k+1))=d_k(\omega_{i_2,j_0-1}(k+1))=c_{j_0-1}$.
Since the last $k$ components of $\omega_{i_1,j_0-1}(k+1)$ and $\omega_{i_2,j_0-1}(k+1)$ are the same, the one-to-one-ness of $d_k$ with respect to the first variable implies that $x_{i_1,j_0-1}=x_{i_2,j_0-1}$. Therefore
$\omega_{i_1,j_0-1}(n+k)=\omega_{i_2,j_0-1}(n+k)$, which
contradicts the minimality of $j_0$. This establishes the theorem.
\end{proof}

In general we see that applying the inverse of a homomorphism to an vertex disjoint cycle in $B_n(q)$ creates multiple cycles in $B_{n+k}(q)$. If $B_n(q)$ is partitioned into vertex disjoint cycles then the inverse homomorphism naturally induces a partition of $B_{n+k}(q)$ into vertex disjoint cycles.

\subsection{The Binary case}\label{S:binaryPropD}
We treat here the binary case separately because its simplicity allows for a more concise characterization of the shape of homomorphisms with property (D).
\begin{theorem}\label{T:propertyD}
A necessary and sufficient condition for a homomorphism $D_{n,k}$ from $B_{n+k}$ to $B_n$
to have property (D) is that

\[d_k(x_1,\ldots,x_{k+1})=x_1+h(x_2,\ldots,x_k)+x_{k+1},
\]
where $h(x_2,\ldots,x_k)$ is any Boolean function of $k-1$ variables.
\end{theorem}
\begin{proof}
By Theorem \ref{T:genpropertyD} we only need to show that a binary function $d_k$ is bijective with respect to the first and last variables if and only if it has the form claimed in this Theorem. In effect, if $d_k(x_1,\ldots,x_{k+1})$ is bijective in $x_1$ and in $x_{k+1}$ then it satisfies the equations
\[
d_k(\bar{x}_1,x_2,\ldots,x_{k+1})=1-d_k(x_1,x_2,\ldots,x_{k+1})=d_k(x_1,x_2,\ldots,\bar{x}_{k+1}).
\]
\noindent So that $d_k(\bar{x}_1,x_2,\ldots,x_k,\bar{x}_{k+1})=d_k(x_1,x_2,\ldots,x_k,x_{k+1})$.

Therefore for each fixed set of values for
$x_2,\ldots,x_k$, $d_k(x_1,\ldots,x_{k+1})=d_{x_2,\ldots,x_k}(x_1,x_{k+1})$ is either
$x_1+x_{k+1}$ or $x_1+x_{k+1}+1$. This can be rephrased to establish the necessity. The converse is obvious because $d_k$ is linear in the first and last variables.
\end{proof}

This elegant form of $d_k$ is mainly due to the ``lack'' of terms in $\mathbb{Z}_2$. While Theorem~\ref{T:genpropertyD} shows that $d_k=\alpha x_1+h(x_2,\ldots,x_k)+\beta x_{k+1}$ is sufficient for property (D), the following example illustrates why property (D) homomorphisms can not be all written in such a simple form even for $q=3$. In fact all the twelve $3\times3$ Latin squares can be written in function form as $f(b_1,b_2)=\alpha_1b_1+\alpha_2b_2+\alpha_3$ where $b_i,\alpha_i\in\mathbb{Z}_3$, $\alpha_1\neq0$ and $\alpha_2\neq0$. Fr the values $0$, $1$, $2$ of $x_2$ let $d_k(x_1,x_2,x_3)$ be respectively $x_1+x_2+x_2$, $2x_1+x_3$ and $x_1+2x_2+2x_3$. Then $d_k$ has property (D) by Theorem~\ref{T:genpropertyD} but it is not linear in either $x_1$ or $x_{k+1}$, despite the simple form of Latin squares. Notice that when $q>3$ most Latin squares are already nonlinear.

While the only binary homomorphism for $k=1$ is Lempel's D-morphism (and
its bitwise complement), there are essentially two homomorphisms
for $k=2$ that are induced by the functions $d^{(1)}=x_1+x_3$ and
$d^{(2)}=x_1+x_2+x_3$. Note that the former is just the D-morphism
iterated twice. The only other two homomorphisms are
bitwise complements of $d^{(1)}$ and $d^{(2)}$. The cases $k\geq3$ allow for nonlinear homomorphisms such as $d(x_1,\ldots,x_4)=x_1+x_2x_3+x_4$.

Let $\mathbf{c}=[c_1\ldots c_l]$ be an arbitrary but fixed cycle in $B_n$, started at a fixed word, say, $0\ldots0$. Then for each
homomorphism $D_{n,k}$ with property (D), $\mathbf{c}$ defines a
map $\mathcal{D}_{\mathbf{c}}$ on the set $\mathbb{Z}_2^k$ as follows. $\mathcal{D}_{\mathbf{c}}(z_1\ldots
z_k)$ is the suffix of length $k$ of
$\mathcal{D}_{\mathbf{c}}^{-1}\mathbf{c}$ started at the string $z_1\ldots z_k$. The inverse
image is generated by the recursive equation
\[z_i=c_{i-k}+z_{i-k}+h(z_{i-k+1},\ldots,z_{i-1});\;i=k+1,\ldots,k+l,
\]
where $h$ is as in Theorem ~\ref{T:propertyD} and
$z_1,\ldots,z_k$ are the required initial conditions. It can be seen that property (D) implies that $\mathcal{D}_{\mathbf{c}}$ is a bijection.
%We state this as
%\begin{proposition}
%Given a homomorphism $D_{n,k}$ with property (D) from $B_{n+k}$ to $B_n$, every vertex disjoint cycle $\mathbf{c}$ in $B_n$ with a fixed initial vertex defines a permutation on the set $\mathbb{Z}_2^k$ where the image of a $k$-word $\mathbf{x}$ is the $k$-word suffix of the inverse image of $\mathbf{c}$ by $D_{n,k}$.
%\end{proposition}
%
When the D-morphism is used, any De~Bruijn cycle $\mathbf{b}_n$ yields the
identity permutation on the set $\{0,1\}$. This is a restatement
of the fact that the inverse image of any De ~Bruijn cycle
$\mathbf{b}_n$ under the $D$-morphism makes two dual cycles in
$B_{n+1}$. Since a binary De ~Bruijn cycle necessarily has an even number of ones, this follows immediately from Fact~ 1 above. The next proposition
concerns the function $d^{(2)}$ defined above.

\begin{proposition}\label{P:span3homomorphism}
For any integer $n\geq1$ and any De ~Bruijn cycle $\mathbf{b}_n=[b_1\ldots b_{2^n}]$, the
homomorphism induced by the Boolean function
$d^{(2)}(x_1,x_2,x_3)=x_1+x_2+x_3$ defines a permutation of the
set of seeds $\{00,01,10,11\}$ with exactly one fixed point
$z_1z_2$ obtained by

$z_1=a_{0}+\delta_{\tilde{n},0}a_{1}+\delta_{\tilde{n},1}a_2$,
$z_2=a_1+\delta_{\tilde{n},0}a_2+\delta_{\tilde{n},1}a_{0}$,

\noindent where $\tilde{n}=n\mod2$,
$a_j=a_j^{(n)}:=\sum_i b_{3i+j}\mod2;j=0,1,2$, the sum is taken
over the range of indices of $\mathbf{b}_n$ ($1\leq3i+j\leq2^n$),
and addition in the index of $a_{j}$ is taken modulo $3$.
\end{proposition}

In other words, exactly one of the four sequences that form the
preimage of $\mathbf{b}_n$ is a closed cycle in $B_{n+2}$. As a result, the other sequences
together form one cycle of length $3\cdot2^{n}$.

\begin{example}
Let $\mathbf{b}_3=[00011101]$. We see that $\tilde{n}=1$, $a_0=1$, $a_1=1$ and $a_2=0$ so that the fixed point is $z_1z_2=10$. Indeed
the inverse image by $d^{(2)}$ gives the following four sequences.

$$\underline{00}000100\underline{01};\;\;\underline{01}101001\underline{11};\;\;\underline{10}110010\underline{10};\;\;
\underline{11}011111\underline{00}.$$

\noindent So the fixed point gives the only cycle of length $8$ while the other
three sequences make the following cycle of length $24$:

$[\underline{00}000100\underline{01}101001\underline{11}011111]$.
\end{example}

\begin{proof} (of Proposition ~\ref{P:span3homomorphism}) Let
$\bar{i}=i\mod 3$. Iterating the relation
$z_i=b_{i-2}+z_{i-1}+z_{i-2}$, which is satisfied by the
sequence $\{z_i\}_{i=3}^{2^n+2}$, we get (for all $i$ in the range
of the latter sequence)
\[ z_i=\sum_{j=0}^{\lfloor
i/3\rfloor-1}(b_{3j+\bar{i}}+b_{3j+\bar{i}+1})+
\delta_{\bar{i},1}z_1+\delta_{\bar{i},2}z_2+\delta_{\bar{i},0}(z_1+z_2).
\]
where we define  $b_0$ to be zero. Note that $2^n\mod3=1$ or $2$ when $n$ is even or odd
respectively. For each of these two cases, using the above
recursive equation and the requirement $z_{2^n+j}=z_j;\;j=1,2$
yields two linearly independent equations whose unique solution is
as claimed.
\end{proof}

Shifting $\mathbf{b}_n$ by a number that is not a
multiple of $3$ permutes the numbers $a_j,\;j=0,1,2$. So it
changes the permutation but still keeps one fixed point. This
result is interesting because it is independent of the De ~Bruijn
cycle used. The permutation induced by $d^{(1)}$ may or may not have a
fixed point, depending on $\mathbf{b}_n$. As a result,
The function $d^{(2)}$ can be used to generate De ~Bruijn cycles
recursively by joining the shorter cycle (the one started at the
fixed point) to the long cycle made of the other three starting
values. Proposition ~\ref{P:span3homomorphism} describes the way
to identify the two starting digits of the short cycle.
There is no simple way to identify a pair of conjugate words to perform this cross-join operation a priori, for example the alternating strings may or may not be on the same cycle. Note that the existence of a word without a conjugate on the shorter cycle is guaranteed because otherwise the cycle must be a De ~Bruijn cycle, see \cite{fred1982}. Consequently, one can find a cross-join pair by only searching the shorter cycle for a word without a conjugate there. This search takes $O(N^2)$ in the worst case, where $N=2^n$ is the length of the short cycle. This is manageable for small to medium word size $n$.


\begin{thebibliography}{99}
%\bibitem{BEST51} van Aardenne-Ehrenfest,T.,  De Bruijn, N. G.: Circuits and
%Trees in Ordered Linear Graphs. Simon Steven 28, 203-217 (1951)
\bibitem{AA08} A. Alhakim, and M. Akinwande, A Multiple Stream Generator Based on De ~Bruijn Digraph Homomorphisms, Journal of Statistical Computing and Simulation, in press.
\bibitem{AA08w} A. Alhakim, \verb"www.clarkson.edu/~aalhakim/code/recursiveDeBruijn.html".
\bibitem{Annex97}  F. S. Annexstein, Generating De Bruijn Sequences: An Efficient
Implementation, IEEE Transactions on Computers 46, 2 (1997); 198-200.
\bibitem{Chan82} A.H. Chan, R.A. Games, and E.L. Key, On the complexities of De~Bruijn sequences, Journal of Combinatorial Theory Series A 33 (1982); 233-246.
\bibitem{Chang99}  T. Chang, B. Park and Y. H. Kim, An Efficient
Implementation of the $D$-Homomorphism for Generation of De~Bruijn Sequences,
IEEE Transactions on Information Theory 45, 4 (1999); 1280-1283.
\bibitem{Chen1995}  C. Chen, and J. Chen, A homomorphism of the De Bruijn graphs and its applications, IEEE first international conference on algorithms and architectures for parallel processing,  (1995); 465-470.
\bibitem{deBruijn46}  N. G. De Bruijn, A Combinatorial Problem, Koninklijke Nederlandse Akademie v. Wetenschappen 49  (1946); 758–764.
\bibitem{EtzionLempel84} T. Etzion and A. Lempel, Construction of De Bruijn sequences of minimal complexity, IEEE Transactions on Information Theory 30, 5 (1984); 1280-1283.
\bibitem{fred1982}  H. Fredricksen,   A Survey of Full Length Nonlinear Shift Register Cycle Algorithms, SIAM Review, 24, 2 (1982); 195-221.
\bibitem{Golomb1967} S. Golomb,  Shift Register Sequences, San Francisco, Holden-Day, 1967.
\bibitem{Knuth2005} D. Knuth,  The Art of Computer Programming,  4, Fascicle 2: Generating All Tuples and Permutations (Art of Computer Programming) online version.

\bibitem{Leach1960}  E. B. Leach, Regular Sequences and Frequency Distributions, Proceedings of the American
Mathematical Society 11(1960); 566-574.
\bibitem{Lempel70}  A. Lempel, On a Homomorphism of the De Bruijn Graph and Its Applications to the Design of Feedback Shift Registers, IEEE Transactions on Computers  C-19, 12 (1970); 1204-1209.

\bibitem{Martin34}  M. H. Martin, A Problem in Arrangements, Bulletin of the American Mathematical Society 40 (1894); 859--864.
%\bibitem{Massey69} J. L. Massey, Shift Register Synthesis and BCH
%Decoding. IEEE Transactions on Information Theory, vol. IT-15 No
%1, pp. 122-127, (1969).
\bibitem{MK92} G. Matsumoto  and Y. Kurita, Twisted GFSR generators,
ACM Transactions on Modeling and computer simulation 2 (1994); 179--194.
%\bibitem{Morena05} Eduardo Moreno, De Bruijn sequences and De Bruijn graphs for a general language
%Information Processing Letters, vol. 96, Issue 6, 214-219, 31 (2005).
%\bibitem{Reed62} I.S. Reed and R. M. Stewart, Note on the Existence of Perfect
%Maps, IRE Transactions on Information Theory, vol. IT-8, pp.
%10-12, (1962).
\bibitem{Ronse84}  C. Ronse, Feedback Shift Registers (Lecture
Notes in Computer Science) Springer-Verlag, Berlin, 1984.
\bibitem{Ruskey96} F. Ruskey, Combinatorial generation, Unpublished manuscript, Working Version (1h) (1996).
%\bibitem{Savage97}  Savage, C.: A Survey of Combinatorial Gray Codes. SIAM Review 39, 4, 605-629 (1997)
\bibitem{Sloane}  N. Sloane, On-Line Encyclopedia of Integer Sequences.
\end{thebibliography}
\end{document}